\input amstex
\magnification =\magstep 1
\documentstyle{amsppt}
\pageheight{9truein}
\pagewidth{6.5truein}
\NoRunningHeads
\baselineskip=16pt

\topmatter

\title Sums of element orders in groups of odd order
\endtitle

\author Marcel Herzog*, Patrizia Longobardi** and Mercede Maj**
\endauthor

\affil *School of Mathematical Sciences \\
       Tel-Aviv University \\
       Ramat-Aviv, Tel-Aviv, Israel
{}\\
       **Dipartimento di Matematica \\
       Universit\`a di Salerno\\
       via Giovanni Paolo II, 132, 84084 Fisciano (Salerno), Italy
\endaffil
\abstract Denote by $G$ a finite group and by $\psi(G)$ the sum of
element orders in $G$. If $t$ is a positive integer, denote by $C_t$
the cyclic group of order $t$ and write $\psi(t)=\psi(C_t)$.
In this paper we proved the following {\bf Theorem A}: Let $G$ be a
non-cyclic group of odd order $n=qm$, where $q$ is the smallest
prime divisor of $n$ and $(m,q)=1$. Then the following statements hold.
{\bf (1)} If $q=3$, then
$\frac {\psi(G)}{\psi(|G|)}\leq \frac {85}{301}$, and equality holds if
and only if $n=3\cdot 7\cdot m_1$ with $(m_1,42)=1$  and
$G=(C_7\rtimes C_3)\times C_{m_1}$, with $C_7\rtimes C_3$ non-abelian.
{\bf (2)} If $q>3$, then
$\frac {\psi(G)}{\psi(|G|)}\leq \frac {p^4+p^3-p^2+1}{p^5+1}$,
where $p$ is the smallest prime bigger than $q$ and  equality holds if
and only if $n=qp^2m_1$ with $(m_1,p!)=1$  and
$G=C_q\times C_p\times C_p \times C_{m_1}$.


MSC[2010] 20D60; 20E34
\endabstract

\thanks This work was supported by the National Group for Algebraic and
Geometric Structures, and their Applications (GNSAGA - INDAM), Italy.
The first author  is grateful to the Department of
Mathematics of the University of Salerno
for its hospitality and
support, while this investigation was carried out.
\endthanks
\endtopmatter

\document
\heading 1. Introduction \\
\endheading
  
The problem of detecting structural properties of periodic groups by looking at
the orders of their elements has been considered by various authors, from many 
different points of view. In [1], H. Amiri, S.M. Jafarian Amiri and I.M. Isaacs
introduced the function
$$\psi(G)=\Sigma \{o(x)\mid x\in G\},$$
where $o(x)$ denotes the order of the element $x$ of the finite group $G$. They
proved the following theorem:
\proclaim {Theorem 1} Let $G$ be a group of order $n$, and denote by $C_n$ the
cyclic group of order $n$. Then $\psi(G)\leq \psi(C_n)$, with equality if and only if
$G=C_n$.
\endproclaim
Thus $C_n$ is the unique group of order $n$ with the largest value of $\psi(G)$ for
groups of that order. Following this publication, in several recent papers 
it was proved 
that certain classes of finite groups $G$ can be characterized using the function 
$\psi(G)$ (see [2],[3],[5],[6],[8],[9],[10],\linebreak
[12],[14],[15] and [18]).

The problem of determining the second largest value 
of $\psi$ for groups of order $n$ and finding
groups satisfying that condition was discussed in several papers, for example [4],[7],
[14]
and [17]. Recently S.M. Jafarian Amiri and M. Amiri in [7] (see also [4]) and R. Shen,
G. Chen and C. Wu in [17] studied finite groups $G$ of order $n$ with the second largest
value of $\psi(G)$, and obtained information about the structure of $G$ if 
$n=p_1^{\alpha_1}\dotsb p_t^{\alpha_t}$, where $p_i$ are primes satisfying 
$p_1<\dotsb
<p_t$, $\alpha_i$ are positive integers and $\alpha_1>1$. 
Our aim is to continue this investigation in the case when $\alpha_1=1$.

In our paper [11] we obtained the following result which improves Theorem 1:
\proclaim {Theorem 2} If $G$ is a non-cyclic group of order $n$, then $\psi(G)\leq
\frac 7{11}\psi(C_n)$. Moreover, this upper bound is best possible. 
\endproclaim
Recently we determined all groups $G$ of order $n$ satisfying $\psi(G)=\frac
7{11}\psi(C_n)$. In fact, we proved in [14] the following theorem:
\proclaim {Theorem 3} Let $G$ be a non-cyclic group of order $n$. Then $\psi(G)=\frac
7{11}\psi(C_n)$ if and only if $n=4m$ with $m$ an odd integer and $G=C_2\times C_2\times
C_m$.
\endproclaim

In the paper [13], we studied groups of order $n=p_1^{\alpha_1}\dotsb p_t^{\alpha_t}$
with $p_i$ and $\alpha_i$ as stated above,
where $p_1=2$ and $\alpha_1=1$. We proved the following two theorems.
\proclaim {Theorem 4} Let $G$ be a non-cyclic group of order $n=2m$, where $m$ is
an odd integer. Then $\psi(G)\leq \frac {13}{21}\psi(C_n)$. Moreover, $\psi(G)=
\frac {13}{21}\psi(C_n)$ if and only if $G=S_3\times C_{n/6}$, where $n=6m_1$
with $(m_1,6)=1$ and $S_3$ is the symmetric group on three letters.
\endproclaim 

Thus the groups $G=S_3\times C_{n/6}$, where $n=6m_1$ and $(m_1,6)=1$, are the groups
with the second largest value of $\psi$ for groups of order $2m$, with $m$ odd.

The second result in [13] is the following theorem.
\proclaim {Theorem 5} Let $\Delta_n$ be the set of non-cyclic groups of the fixed
order $n=2m$, where $m$ is an odd integer, and suppose that $m=p_1^{\alpha_1}\dotsb
p_t^{\alpha_t}$, where $p_i$ are distint primes and $\alpha_i$ are positive integers
for all $i$. If $G\in \Delta_n$, then 
$$\psi(G)\leq (\frac 13+\frac {2l}{3\psi(l)})\psi(C_n),$$
where $l=\min(p_i^{\alpha_i}\mid i\in \{1,\dots,t\})$ and $\psi(l)=\psi(C_l)$.

Moreover, $G\in \Delta_n$ satisfies 
$\psi(G)=(\frac 13+\frac {2l}{3\psi(l)})\psi(C_n)$ if and only if 
$G=D_{2l}\times C_{n/2l}$, where $D_{2l} $ denotes the dihedral group of order $2l$.
\endproclaim

Thus if $n=2m$, with $m$ odd, then $D_{2l}\times C_{n/2l}$ is the unique group of 
order $n$ with the second largest value of $\psi$ for groups of that order. In
particular, if $l=3$, then $D_{2l}=S_3$ 
and $\frac 13+\frac {2l}{3\psi(l)}=\frac 13+\frac 6{3\cdot 7}=
\frac {13}{21}$, in agreement with Theorem 4.
 
Using Theorem 2, Theorem 3, Theorem 4 and results from [7] and [17], we could
obtain the following result:

\proclaim {Theorem 6} Let $G$ be a non-cyclic group of even order 
$n=2^{\alpha}m$, with $m$ an odd integer. Then the following statements hold.
\roster
\item If $\alpha =1$, then $\frac {\psi(G)}{\psi(C_n)}\leq \frac {13}{21}=
\frac {\psi(S_3)}{\psi(C_6)}$,
\item If $\alpha =2$, then $\frac {\psi(G)}{\psi(C_n)}\leq \frac 7{11}=
\frac {\psi(C_2\times C_2)}{\psi(C_4)}$,
\item If $\alpha =3$, then $\frac {\psi(G)}{\psi(C_n)}\leq \frac {27}{43}=
\frac {\psi(Q_8)}{\psi(C_8)}$,
\item If $\alpha \geq 4$, then $\frac {\psi(G)}{\psi(C_n)}\leq 
\frac {2^{2\alpha +3}+7}{7(2^{2\alpha +1}+1)}=
\frac {\psi(G_1)}{\psi(C_{2^{\alpha}\cdot 3})}$,
\endroster
where $Q_8$ is the quaternion group of order $8$ and $G_1=\langle a,b\mid 
a^{2^\alpha}=b^3=1, b^a=b^{-1}\rangle $.

In particular, these upper bounds are best possible.
\endproclaim

In the paper [14], we started the investigation in the odd case. By Theorem 3
in [11] we have
\proclaim {Theorem 7} If $G$ is a non-cyclic group of order $n$ and if $q$ 
denotes the smallest prime divisor of $n$, then 
$\frac {\psi(G)}{\psi(C_n)}<\frac 1{q-1}$.
\endproclaim
 
Theorem 7 implies the following corollary.

\proclaim {Corollary 8} If $G$ is a non-cyclic group of odd order $n$, then 
$\frac {\psi(G)}{\psi(C_n)}<\frac 12$.
\endproclaim

In [14] we proved the following theorem, which  yields a better result
in the odd case.

\proclaim {Theorem 9}  If $G$ is a non-cyclic group of order $n$ and if $q$
denotes the smallest prime divisor of $n$, then 
$$\frac {\psi(G)}{\psi(C_n)}\leq \frac {q^4+q^3-q^2+1}{q^5+1}$$
and the equality holds if and only if $n=q^2k$ with $(k,q!)=1$ and 
$G=C_q\times C_q\times C_k$.
\endproclaim

As shown in Lemma 2.1 of Section 2, the 
function $f(x)= \frac {x^4+x^3-x^2+1}{x^5+1}$ is strongly decreasing for $x\geq 2$.
Notice that $f(2)=\frac 7{11}$ and $f(3)=\frac {25}{61}$.
Thus Theorem 9 implies the following result.

\proclaim {Theorem 10} If $G$ is a non-cyclic group of order $n$,
then the following statements hold.
\roster
\item If $G$ is of even order, then   $\frac {\psi(G)}{\psi(C_n)}\leq \frac 7{11}$
and the equality holds if and only if $n=4k$ with $(k,2)=1$ and 
$G=C_2\times C_2\times C_k$.
\item If $G$ is of odd order, then   $\frac {\psi(G)}{\psi(C_n)}\leq \frac {25}{61}$
and the equality holds if and only if $n=9k$ with $(k,6)=1$ and
$G=C_3\times C_3\times C_k$.
\endroster
\endproclaim

\demo{Proof} Let $q$ denote the smallest prime divisor of $n$. If $n$ is even, then
$q=2$ and by Theorem 9 $\frac {\psi(G)}{\psi(C_n)}\leq f(2)$, with equality as 
required. If $n$ is odd, then $q\geq 3$ and by Theorem 9 
$\frac {\psi(G)}{\psi(C_n)}\leq f(q)$. It follows by Lemma 2.1 that $f(q)\leq f(3)$,
so  $\frac {\psi(G)}{\psi(C_n)}\leq f(3)=\frac {25}{61}$, as required. The condition
for equality is again as required.
\qed
\enddemo

Hence if $G$ is a non-cyclic group of odd order, then 
$\frac {\psi(G)}{\psi(C_n)}\leq \frac {25}{61}<\frac 12$ and the upper bound 
$\frac {25}{61}$ is best possible.
  
In this paper we are looking for a result similar to Theorem 5 for groups of odd order
$n$. Denote by $q$ the smallest prime divisor of $n$. As in the case when $q=2$,
we start with groups of order $n=qm$, where $q>2$ and $(q,m)=1$. We denote by
$C_7\rtimes C_3$ the non-abelian group of order $21$. We prove the following theorem:

\proclaim {Theorem A}  Let $G$ be a
non-cyclic group of odd order $n=qm$, where $q$ is the smallest
prime divisor of $n$ and $(m,q)=1$. Then the following statements hold.
\roster
\item If $q=3$, then
$\frac {\psi(G)}{\psi(C_n)}\leq \frac {85}{301}$, and equality holds if
and only if $n=3\cdot 7\cdot m_1$ with $(m_1,42)=1$  and
$G=(C_7\rtimes C_3)\times C_{m_1}$, with $C_7\rtimes C_3$ non-abelian.
\item  If $q>3$, then
$\frac {\psi(G)}{\psi(C_n)}\leq \frac {p^4+p^3-p^2+1}{p^5+1}$,
where $p$ is the smallest prime bigger than $q$ and  equality holds if
and only if $n=qp^2m_1$ with $(m_1,p!)=1$  and
$G=C_q\times C_p\times C_p \times C_{m_1}$.
\endroster
\endproclaim

Theorem A implies the following corollary.
\proclaim {Corollary B}  Let $G$ be a non-cyclic group of odd order $n=qm$,
where  $q$ is the smallest
prime divisor of $n$ and $(m,q)=1$. Then the following statements hold.
\roster
\item If $q = 3$, then
$$\frac {\psi(G)}{\psi(C_n)}\leq \frac {85}{301}\approx 0.282,$$
and  equality holds if
and only if $q=3$, $n=3\cdot 7\cdot m_1$ with $(m_1,42)=1$  and
$G=(C_7\rtimes C_3)\times C_{m_1}$.
\item If $q>3$, then
$$\frac {\psi(G)}{\psi(C_n)}\leq \frac {337}{2101}\approx 0.160,$$
and  equality holds if
and only if $q=5$, $n=5\cdot 7^2\cdot m_1$ with $(m_1,7!)=1$  and
$G=C_5\times C_7\times C_7 \times C_{m_1}$.
\endroster
\endproclaim
 
Finally, Theorem A and results of [7], [14] and [17] imply the following corollary.
\proclaim {Corollary C} Let $G$ be a non-cyclic group of odd order $n=mq^{\alpha}$,
where $q$ is the smallest prime divisor of $n$, $\alpha$ is a positive integer and
$(m,q)=1$. Moreover, let $p$ be the smallest prime bigger than $q$. Then the following
statements hold.
\roster
\item If $\alpha =1$ and $q=3$, then  
$\frac {\psi(G)}{\psi(C_n)}\leq \frac {85}{301}=\frac {\psi(C_7\rtimes C_3)}
{\psi(C_{21})}$.
\item If $\alpha =1$ and $q\geq 3$, then
$\frac {\psi(G)}{\psi(C_n)}\leq \frac {p^4+p^3-p^2+1}{p^5+1}
=\frac {\psi(C_p\times C_p)}{\psi(C_{p^2})}$.
\item If $\alpha =2$, then
$\frac {\psi(G)}{\psi(C_n)}\leq \frac {q^4+q^3-q^2+1}{q^5+1}
=\frac {\psi(C_q\times C_q)}{\psi(C_{q^2})}$.
\item  If $\alpha =3$, then
$\frac {\psi(G)}{\psi(C_n)}\leq \frac {q^6+q^3-q^2+1}{q^7+1}
=\frac {\psi(C_q\times C_{q^2})}{\psi(C_{q^3})}=\frac {\psi(M(q^3))}{\psi(C_{q^3})}$,
\newline
where $M(q^3)=\langle a,b\mid a^{q^2}=b^q=1,\ a^b=a^{q+1}\rangle$.
\item  If $\alpha \geq 4$, then
$\frac {\psi(G)}{\psi(C_n)}\leq \frac {q^{2\alpha}+q^3-q^2+1}{q^{2\alpha+1}+1}
=\frac {\psi(C_q\times C_{q^{\alpha -1}})}{\psi(C_{q^{\alpha}})}$.
\endroster
These upper bounds are best possible.
\endproclaim

The proof of Theorem A will be achieved in two steps. The first step is
\proclaim {Proposition D} Let
$G$ be a non-cyclic group of odd order $n=qm$, where  $q$ is the smallest
prime divisor of $n$ and $(m,q)=1$. Then the following statements hold.
\roster
\item If $q=3$, then
$$\frac {\psi(G)}{\psi(C_n)}\leq \frac {85}{301}.$$
\item If $q>3$, then
$$\frac {\psi(G)}{\psi(C_n)}\leq\frac {p^4+p^3-p^2+1}{p^5+1},$$
where $p$ denotes the smallest prime bigger than $q$.
\endroster
\endproclaim
The second step is
\proclaim {Proposition E} Let
$|G|=n=qm$, where $q$ is the smallest
prime divisor of $n$ and $(m,q)=1$. Then the following statements hold.
\roster
\item If $q=3$, then
$$\frac {\psi(G)}{\psi(C_n)}=\frac {85}{301}$$
if
and only if $n=3\cdot 7\cdot m_1$ with $(m_1,42)=1$  and
$G=(C_7\rtimes C_3)\times C_{m_1}$.
\item If $q>3$ and if $p$ denotes the smallest prime bigger than
$q$, then
$$\frac {\psi(G)}{\psi(C_n)}=\frac {p^4+p^3-p^2+1}{p^5+1}$$
if and only if $n=qp^2m_1$ with $(m_1,p!)=1$ and
$G=C_q\times C_p\times C_p \times C_{m_1}$.
\endroster
\endproclaim

It is clear that Theorem A follows from Propositions D and E. These
propositions will be proved in Sections 3 and 4, respectively. Section
2 will be devoted to preliminary results and Corollary B will be proved in
Section 5.

Our notation will be the usual one (see for example [16]). In particular, 
$\varphi(n)$ denotes the Euler's function on $n$. If $k$ is a positive integer,
we write $\psi(k)=\psi(C_k)$, and, following [17], $\lambda(k)=\frac {\psi(k)}k$.
Finally, if $X$ is a subset of $G$, we shall define $\psi(X)=\Sigma_{x\in X}o(x)$.

\heading 2. Preliminary results \\
\endheading

First we define two functions of real variables which will play important
roles in our proofs.
\definition{Definition} Let $q>1$ be a positive integer. Then
$$f(x)=\frac {x^4+x^3-x^2+1}{x^5+1}\quad \text{and}\quad
g_q(x)=\frac {x^2-x+1+x(q^2-q)}{(x^2-x+1)(q^2-q+1)}.$$
\enddefinition

Our first result is the following lemma.
\proclaim{Lemma 2.1} The functions $f$ and $g_q$ are strongly decreasing for
$x\geq 2$.
\endproclaim
\demo{Proof} We have $g_q'(x)=\frac
{(q^2-q)(1-x^2)}{(q^2-q+1)(x^2-x+1)^2}
<0$ for $x\geq 2$. Furthermore, if $x \geq 2$, then
$f(x)=\frac {x^3-x+1}{x^4-x^3+x^2-x+1}$ and
$f'(x)=\frac {-xh(x)}{(x^4-x^3+x^2-x+1)^2}$,
where $h(x)=x^5-4x^3+8x^2-7x+2>0$. Hence also $f'(x)<0$ for $x\geq 2$.
\qed
\enddemo

In sequel, we shall often use without reference the fact that if $m$ and $n$ 
are integers
satisfying $m>n\geq 2$, then $f(m)<f(n)$ and $g_q(m)<g_q(n)$.

Our second result is the
following important proposition.
\proclaim{Proposition 2.2} Let $q\geq 3$ be a prime. Denote by $p$
the smallest prime bigger than $q$ and by
$q_1$ the smallest
prime congruent to $1 \ (mod\ q)$.
Then the following statements hold.
\roster
\item  If $q=3$, then $q_1=7$, $p=5$ and
$$\frac {85}{301}=g_3(q_1)>f(p)=\frac {121}{521}.$$
\item If $q>3$, then
$$g_q(q_1)=\frac {q_1^2-q_1+1+q_1(q^2-q)}{(q_1^2-q_1+1)(q^2-q+1)}<f(p)=\frac
{p^4+p^3-p^2+1}{p^5+1}.$$
\endroster
\endproclaim
\demo{Proof} 
{\bf (1)} If $q=3$, then  $p=5$ and $q_1=7$. Thus
$g_3(q_1)=\frac {7^2-7+1+7(3^2-3)}{(7^2-7+1)(3^2-3+1)}=\frac {85}{301}$ and
$f(p)=\frac {5^4+5^3-5^2+1}{5^5+1}=\frac {121}{521}<\frac {85}{301}$, as
required.

{\bf (2)} By our assumption $q>3$.
Clearly $q_1\geq 2q+1$ and by the Bertrand's conjecture
$p\leq 2q-3$. As shown in Lemma 2.1, the functions $g_q(x)$
and $f(x)$ are both decreasing for $x\geq 2$, which implies that
$$g_q(q_1)\leq g_q(2q+1)\quad \text{and}\quad f(p)\geq f(2q-3).$$
Therefore it suffices to prove that $g_q(2q+1)<f(2q-3)$. So we need to prove that
$$\frac {(2q+1)^2-(2q+1)+1+(2q+1)(q^2-q)}{((2q+1)^2-(2q+1)+1)(q^2-q+1)}<
\frac {(2q-3)^4+(2q-3)^3-(2q-3)^2+1}{(2q-3)^5+1}.$$
Since $\frac {a+1}{b+1}>\frac ab$ if $b>a$, it suffices to prove that
$$\frac {2q^3+3q^2+q+1}{4q^4-2q^3+3q^2+q+1}<\frac {(2q-3)^2+(2q-3)-1}{(2q-3)^3}$$
or
$$\frac {2q^3+3q^2+q+1}{4q^4-2q^3+3q^2+q+1}<\frac {4q^2 -10q
+5}{8q^3-36q^2+54q-27}.$$
This inequality is equivalent to the inequality
$$16q^6-48q^5+8q^4+80q^3-63q^2+27q-27<16q^6-48q^5+52q^4-36q^3+9q^2-5q+5,$$
or
$$116q^3+32q<44q^4+72q^2+32,$$
which is certainly true if $q>3$.

The proof of the proposition is now complete.
\qed
\enddemo
In the next lemma, we collect some basic information
about the function $\psi$.
\proclaim {Lemma 2.3}  Let $G$ denote a finite group, $p,p_i$ denote primes
and $n,m,\alpha_i$ denote positive integers. Then the following statements hold.
\roster
\item"({\bf 1})" {\rm ([11], Lemma 2.9(1))} $\psi(C_{p^m})=\frac {p^{2m+1}+1}{p+1}=
\frac {p(p^{m})^2+1}{p+1}$.
\item"({\bf 2})" {\rm ([11], Lemma 2.2(3))} If $G=A\times B$, where $A,B$
are subgroups of $G$ satisfying
$(|A|,|B|)=1$, then
$\psi(G)=\psi(A)\psi(B)$.
\item"({\bf 3})" {\rm ([11], Lemma 2.9(2))}
If $n=\prod _{i=1}^{i=m}p_i^{\alpha_i}$,
where $p_i\neq p_j$ for
$i\neq j$, then
$\psi(C_n)=\prod _{i=1}^{i=m}\psi(C_{p_i^{\alpha_i}})$.
\item"({\bf 4})" {\rm ([11], Proof of Lemma 2.9(2))} If $q$  and $p$ are the
smallest and the largest prime divisors
of $n$, then
$\psi(C_n)\geq \frac q{p+1}n^2$.
\item"({\bf 5})" {\rm ([1], Corollary B)} If $P$ is a cyclic normal Sylow
subgroup of $G$, then
$\psi(G)\leq \psi(P)\psi(G/P)$, with equality if and only if $P$ is
central in $G$.
\item"({\bf 6})" {\rm ([11], Lemma 2.2(5))} Let $G=P\rtimes F$, where $P$ is
a cyclic $p$-group and $F$ is a group satisfying $|F|>1$ and
$(p,|F|)=1$. Let $Z=C_F(P)$.
Then
$$ \psi(G)=\psi(P)\psi(Z)+|P|\psi(F\setminus Z).$$
\item"({\bf 7})" {\rm ([11], Lemma 2.2(1))} Let $G=P\rtimes F$, where $P$ and
$F$ are as described above.
Then each element of $F$ acts on
$P$ either trivially or fixed-point-freely.
\endroster
\endproclaim

We conclude this section with the following two lemmas which follow easily from 
previous results in [13].

\proclaim{Lemma 2.4} Let $G=(\langle a\rangle \times \langle b\rangle)\rtimes
\langle y\rangle $, where $p$ is an odd prime number and $o(a)=p^{\alpha -1}$ for some
integer $\alpha >1$, $o(b)=p$, $(o(y),p)=1$, $a^y=a$ and $b^y=b^r$ for some integer $r$ 
not congruent to $1$ modulo $p$. Then
$$\psi(G)< \psi(\langle a\rangle \times \langle b\rangle)\psi(\langle y\rangle).$$
\endproclaim

\demo{Proof} By Theorem 3.2 of [13] we have $\psi(G)=(p-1)^2\psi(Z)+p\psi(\langle
a\rangle)\psi(\langle y\rangle)$, where $Z=C_{\langle y\rangle}(\langle b\rangle)$.
On the other hand, 
$$\psi(\langle a\rangle \times \langle b\rangle)=\psi(\langle b\rangle)
+p(\psi(\langle a\rangle)-1)= p^2-p+1-p+p\psi(\langle a\rangle)=(p-1)^2
+p\psi(\langle a\rangle).$$

Thus $\psi(\langle a\rangle \times \langle b\rangle)\psi(\langle y\rangle)=
((p-1)^2+p\psi(\langle a\rangle))\psi(\langle y\rangle)$ and 
$Z\subset \langle y\rangle $ implies that $\psi(Z)<\psi(\langle y \rangle)$.
Hence 
$$\psi(G)=(p-1)^2\psi(Z)+p\psi(\langle a\rangle)\psi(\langle y\rangle)
<((p-1)^2+p\psi(\langle a\rangle))\psi(\langle y\rangle)= 
\psi(\langle a\rangle \times \langle b\rangle)\psi(\langle y\rangle),$$
as required. 
\qed
\enddemo
\proclaim{Lemma 2.5} Let $G$ be a non-cyclic group and $G=P\rtimes \langle y\rangle$,
where $P$ is a non-trivial cyclic $p$-subgroup of $G$ for some prime $p$ and 
$(|P|,|\langle y\rangle|)=1$. If $q$ is the smallest prime divisor of
$|\langle y\rangle|$, then
$$\frac {\psi(G)}{\psi(|G|)}\leq \frac {p^2-p+1+p(q^2-q)}{(p^2-p+1)(q^2-q+1)}.$$
\endproclaim
\demo{Proof} Write $F= \langle y\rangle$ and let $Z=C_F(P)$. Then, by Lemma 2.3(6),
we have $\psi(G)=\psi(P)\psi(Z)+|P|(\psi(F)-\psi(Z))$. Thus
$\frac {\psi(G)}{\psi(|G|)}=\frac {\psi(G)}{\psi(P)\psi(F)}=\frac {\psi(Z)}{\psi(F)}
+\frac {|P|}{\psi(|P|)}(1-\frac {\psi(Z)}{\psi(F)})$.

Write $|P|=p^\alpha$, with $\alpha $ a positive integer, and if $k$ is a positive
integer, let $\lambda(k)=\frac {\psi(k)}k$. By Lemma 2.8 of [13] we have 
$\lambda(p)\leq \lambda(p^\alpha)$, so $\frac 1{ \lambda(p^\alpha)}\leq 
\frac 1{ \lambda(p)}$. Thus we have
$$\frac {\psi(G}{\psi(|G|)}=\frac {\psi(Z)}{\psi(F)}+\frac 1{ \lambda(p^{\alpha})}
\big(1-\frac {\psi(Z)}{\psi(F)}\big)
\leq \frac {\psi(Z)}{\psi(F)}+\frac 1{ \lambda(p)}
\big(1-\frac {\psi(Z)}{\psi(F)}\big)=\frac {\psi(Z)}{\psi(F)}
\big(1-\frac 1{ \lambda(p)}\big)
+\frac 1{ \lambda(p)}.$$

Arguing as in the proof of Proposition 6 of [14], we have
$\frac {\psi(Z)}{\psi(F)}\leq \frac 1{q^2-q+1}$. Hence
$$
\align
\frac {\psi(G)}{\psi(|G|)} &\leq \frac 1{q^2-q+1}\big(1-\frac 1{\lambda(p)}\big)
+\frac 1{\lambda(p)}=\frac 1{q^2-q+1}\big(1-\frac p{p^2-p+1}\big)+\frac p{p^2-p+1}\\
&=\frac {(p^2-2p+1)+p(q^2-q+1)}{(q^2-q+1)(p^2-p+1)}=\frac {p^2-p+1 +p(q^2-q)}
{(q^2-q+1)(p^2-p+1)},
\endalign$$
as required.
\qed
\enddemo

\heading 3. Proof of Proposition D
\endheading

\demo {Proof} We assume that $G$ is a non-cyclic group of
odd order $|G|=n=qm$ with $(m,q)=1$. This implies 
that $G=M\rtimes C$, where $C\simeq C_q$ and $M$ is its normal
complement.

First assume that $q=3$. Then $|G|=3m$ with $(3,m)=1$ and 
suppose that $\frac {\psi(G)}{\psi(|G|)}>\frac {85}{301}$.
Our aim is to reach a contradiction. We argue by induction on $|G|$.

If $G=M\times C$, then
$\psi(G)=\psi(M)\psi (C)$, where $M$ is non-cyclic and $3$ does not divide
$|M|$.
By Theorem 4 of [14] and Lemma 2.1 we have
$\psi(M)\leq \frac {5^4+5^3-5^2+1}{5^5+1}\psi(|M|)=\frac {121}{521}\psi(|M|)$. Hence
$\psi(G)\leq \frac {121}{521}\psi(|G|)<\frac
{85}{301}\psi(|G|)$, a contradiction. Therefore $C$ is not normal in $G$.

Let $p_2$ be the largest prime divisor of $n$. Clearly $p_2\geq 5$. 
By Lemma 2.3(4), we have
$\psi(|G|)\geq \frac 3{p_2+1}n^2$,
so $\psi(G)>\frac {85}{301}\psi(|G|)\geq \frac{255}{301(p_2+1)}n^2$.
Hence there exists $x\in G$ such that
$o(x)>\frac{255}{301(p_2+1)}n$, which implies that
$[G:\langle x\rangle]< \frac {301}{255}(p_2+1)<2p_2$.

Suppose, first, that $p_2$ does not divide $[G:\langle x\rangle]$. Then $\langle
x\rangle$
contains a cyclic Sylow $p_2$-subgroup $P$ of $G$. 
Since $\langle x\rangle\leq
N_G(P)$, it
follows that $[G:N_G(P)]<2p_2$. But $[G:N_G(P)]=1+kp_2$ for some non-negative
integer $k$,
so $1+kp_2<2p_2$ and since $G$ is of odd order, $[G:N_G(P)]=1+p_2$ is
impossible.
Hence
$[G:N_G(P)]=1$ and $P$ is a normal cyclic Sylow $p_2$-subgroup of $G$.
Now Lemma
2.3(5) and our assumptions
imply that $\psi(P)\psi(G/P)\geq \psi(G)>\frac
{85}{301}\psi(P)\psi(|G/P|)$.
Thus $\psi(G/P)>\frac {85}{301}\psi(|G/P|)$. Since $3$ divides $|G/P|$,
it follows by the induction
hypothesis
that $G/P$ is cyclic.
Write $G=P\rtimes F$, where $F$ is isomorphic to
$G/P$.
Then $F$ is cyclic and $3$ divides $|F|$. Since $C$ is not normal in $G$, $3$ does
not divide
$C_F(P)$ and since $P$ is cyclic, it follows by Lemma 2.3(7)
that $p_2$ is congruent to $1$ (mod
$3$). Hence
$p_2\geq 7$ and $g_3(p_2)\leq g_3(7)$. Now
Lemma 2.5 implies that
$$\frac {\psi(G)}{\psi(|G|)}\leq \frac {p_2^2-p_2+1+6p_2}{7(p_2^2-p_2+1)}=g_3(p_2)
\leq  g_3(7)= \frac {85}{301},$$ a contradiction.

Now  suppose that $p_2$ does divide $[G:\langle x\rangle]$. Then $[G:\langle
x\rangle]<2p_2$
implies that $[G:\langle x\rangle]=p_2$. By Theorem 3.1 of [13], $G=P\rtimes F$,
where
$P$ is a normal Sylow $p_2$-subgroup of $G$,
$F$ is a cyclic subgroup of $G$, and either
$P$ is cyclic, or $G$ is nilpotent, or
$G=(\langle a\rangle \times \langle b\rangle)\rtimes \langle y\rangle$, where
$o(a)=p_2^{\alpha -1} $ for some integer $\alpha >1$, $o(b)=p_2$, $(o(y),p_2)=1$,
$a^y=a$ and $b^y=b^r$
for some integer $r$ not congruent to $1$ (mod $p_2$). If $P$ is cyclic, then we get
a contradiction arguing as in the previous paragraph. If $G$ is nilpotent,
then $C$ is normal in $G$,
a contradiction. Finally, in the last case, it follows by Lemma 2.4 that
$\psi(G)<\psi(\langle a\rangle \times \langle b\rangle)\psi(\langle y\rangle)$.
Applying
Theorem 9 to $\langle a\rangle \times \langle b\rangle$
we obtain $\psi(G)<
\frac{p_2^4+p_2^3-p_2^2+1}{p_2^5+1}\psi(|\langle a\rangle \times
\langle b\rangle|)\psi(\langle y\rangle)$. Now
$\psi(|\langle a\rangle \times \langle b\rangle|)\psi(\langle y\rangle)=\psi(|G|)$
and since $p_2\geq 5$, Lemma 2.1 implies that
$\frac{p_2^4+p_2^3-p_2^2+1}{p_2^5+1}\leq \frac {5^4+5^3-5^2+1}{5^5+1}=\frac
{121}{521}$. Thus 
$\psi(G)<\frac {121}{521}\psi(|G|)<\frac {85}{301}\psi(|G|)$, a contradiction.
The proof in the case when $q=3$ is now complete.

We assume, now, that $q>3$.
Then  $|G|=qm$ with $(q,m)=1$ and
suppose that $\frac {\psi(G)}{\psi(|G|)}>\frac {p^4+p^3-p^2+1}{p^5+1}$, where
$p$ is the smallest prime larger than $q$. Denote by $q_1$ the smallest prime
congruent to $1$ modulo $q$. Our aim is to reach a contradiction.
We argue by induction on $|G|$.

If $G=M\times C$, then $\psi(G)=\psi(M)\psi (C)$, where $M$ is non-cyclic
and $q$ does not divide $|M|$. Therefore
$\frac {\psi(G)}{\psi(|G|)}=\frac {\psi(M)}{\psi(|M|)}\frac {\psi(C)}{\psi(|C|)}=
\frac {\psi(M)}{\psi(|M|)}$ and by Theorem 9
$\frac {\psi(M)}{\psi(|M|)}\leq \frac {q_2^4+q_2^3-q_2^2+1}{q_2^5+1}=f(q_2)$,
where $q_2$ is the smallest prime divisor of $|M|$. Obviously $q_2\geq  p$, so
$f(q_2)\leq f(p)$ and
$\frac {\psi(G)}{\psi(|G|)}=\frac {\psi(M)}{\psi(|M|)}\leq \frac
{p^4+p^3-p^2+1}{p^5+1}$, a contradiction. Therefore $C$ is not normal in $G$.

Let $p_2$ be the largest prime divisor of $n$. Then, by Lemma 2.3(4), we have
$\psi(|G|)\geq \frac q{p_2+1}n^2>\frac {q-1}{p_2}n^2$,
so $\psi(G)>\frac {p^4+p^3-p^2+1}{p^5+1}\psi(|G|)> \frac
{(p^4+p^3-p^2+1)(q-1)}{(p^5+1)p_2}n^2$. Hence there exists $x\in G$ such that
$o(x)>\frac{(p^4+p^3-p^2+1)(q-1)}{(p^5+1)p_2}n$, which implies that
$[G:\langle x\rangle]< \frac {(p^5+1)p_2}{(p^4+p^3-p^2+1)(q-1)}$. By Bertrand's
conjecture $p\leq 2q-2$, so $q-1\geq \frac p2$ and
$$[G:\langle x\rangle]< \frac {2(p^5+1)p_2}{(p^4+p^3-p^2+1)p}< 2p_2.$$

First, suppose that $p_2$ does not divide $[G:\langle x\rangle]$. Then $\langle
x\rangle$
contains a cyclic Sylow $p_2$-subgroup $P$ of $G$
and it follows, as in the $q=3$ case, that $P$ is a normal subgroup of $G$. 
Now Lemma
2.3(5) and our assumptions
imply that $\psi(P)\psi(G/P)\geq \psi(G)>\frac
{p^4+p^3-p^2+1}{p^5+1}\psi(P)\psi(|G/P|)$.
Thus $\psi(G/P)>\frac {p^4+p^3-p^2+1}{p^5+1}\psi(|G/P|)$. Since
$q$ divides $|G/P|$, it follows by the inductive 
hypothesis
that $G/P$ is cyclic. 
Write $G=P\rtimes F$, where $F$ is isomorphic to $G/P$. Then $F$ is 
a cyclic group 
and $q$ divides $|F|$.
As shown in the $q=3$ case, $p_2$ is congruent to $1$ (mod
$q$), implying that $p_2\geq q_1$ and $g_q(p_2)\leq g_q(q_1)$.
Since by Proposition 2.2(1) $g_q(q_1)<f(p)$, applying Lemma 2.5 we
obtain
$$\frac {\psi(G)}{\psi(|G|)}\leq \frac
{p_2^2-p_2+1+p_2(q^2-q)}{(p_2^2-p_2+1)(q^2-q+1)}
=g_q(p_2)\leq g_q(q_1)<f(p),$$
a contradiction.

Now suppose that $p_2$ does divide  $[G:\langle x\rangle]$. Since $[G:\langle
x\rangle]<2p_2$,
it follows that $[G:\langle x\rangle]=p_2$. Then, by Theorem 3.1 of [13],
$G=P\rtimes F$, where
$P$ is a normal Sylow $p_2$ subgroup of $G$, $F$ is a cyclic subgroup of $G$, and
either
$P$ is cyclic, or $G$ is nilpotent, or
$G=(\langle a\rangle \times \langle b\rangle)\rtimes \langle y\rangle$, where
$o(a) = p_2^{\alpha -1} $ for some integer $\alpha >1$, $o(b)=p_2$, $(o(y),p_2)=1$,
$a^y=a$ and $b^y=b^r$
for some integer $r$ not congruent to $1$ (mod $p_2$). If $P$ is cyclic, then we get
a contradiction arguing as in the previous paragraph. If $G$ is nilpotent, then $C$ is
normal in $G$,
a contradiction. Finally, in the last case, it follows by Lemma 2.4 that
$\psi(G)<\psi(\langle a\rangle \times \langle b\rangle)\psi(\langle y\rangle)$.
Applying 
Theorem 9 to $\langle a\rangle \times \langle b\rangle$
we obtain $\psi(G)<
\frac{p_2^4+p_2^3-p_2^2+1}{p_2^5+1}\psi(|\langle a\rangle \times
\langle b\rangle|)\psi(\langle y\rangle)$. Now
$\psi(|\langle a\rangle \times \langle b\rangle|)\psi(\langle y\rangle)=\psi(|G))$
and $p_2\geq p$ implies that
$\frac{p_2^4+p_2^3-p_2^2+1}{p_2^5+1}\leq \frac {p^4+p^3-p^2+1}{p^5+1}$.  It follows
that
$\psi(G)<\frac {p^4+p^3-p^2+1}{p^5+1}\psi(|G|)$, a contradiction.

The proof of Proposition D is now complete.               
\enddemo

\heading 3. Proof of Proposition E
\endheading

\demo {Proof} Let $q_1$ be the smallest prime congruent to $1$ (mod $q$).
Recall that $p$ is the smallest prime larger than $q$, 
$f(x)=\frac {x^4+x^3-x^2+1}{x^5+1}$ and 
$g_q(x)=\frac {x^2-x+1+x(q^2-q)}{(x^2-x+1)(q^2-q+1)}$.

First assume that $q=3$.

Then $q_1=7$ and $p=5$. First we show that if $G=C_k\times(C_7\rtimes C_3)$
where $(k,42)=1$, then $\frac {\psi(G)}{\psi(|G|)}=g_3(7)=
\frac {85}{301}$. Indeed,
$\psi(G)=\psi(C_k)\psi(C_7\rtimes C_3)=\psi(C_k)(7^2-7+1+14\cdot 3)=\psi(C_k)85$,
so $\frac {\psi(G)}{\psi(|G|)}=\frac {\psi(C_k)85}{\psi(C_k)\psi(C_7)\psi(C_3)}=
\frac {85}{301}$, as claimed.

Now we turn to the "only if" part.
So suppose that $q=3$, $|G|=n=3m$, $m$ is odd, $(m,3)=1$ and 
$\frac {\psi(G)}{\psi(|G|)}=g_3(7)=\frac {85}{301}$.
Then $G$ is non-cyclic and $G=M\rtimes C$, where $C\simeq C_3$ 
and $3$ does not divide $|M|$. In particular, 
the least prime divisor $r$ of $|M|$ satisfies $r\geq 5$.
Our aim is to show that $G=C_k\times(C_7\rtimes C_3)$
where $(k,42)=1$. We argue by induction on $|G|$.
 
If $G=M\times C$, then $M$ is non-cyclic and $\psi(G)=\psi(M)\psi(C)$. 
Moreover, since $r\geq 5$, Theorem 9 and Lemma 2.1 imply that 
$\psi(M)\leq \frac {(5^3-5+1)(5+1)}{5^5+1}\psi(|M|)
=\frac {121}{521}\psi(|M|)$.
Hence $\frac {\psi(G)}{\psi(|G|)}=\frac {\psi(M)\psi(C)}{\psi(|M|)\psi(C)}
\leq \frac {121}{521}<\frac {85}{301}$, a contradiction.
Therefore $C$ is not normal in $G$.

Let $p_2$ be the largest prime divisor of $|G|$ and let $P$ be a Sylow 
$p_2$-subgroup of $G$. Then, by Lemma 2.3(4), we have 
$\psi(|G|)\geq \frac 3{p_2+1}n^2$ and our assumptions imply that 
$\psi(G)=\frac {85}{301}\psi(|G|)\geq \frac{255}{301(p_2+1)}n^2$.
Hence there exists $x\in G$ such that
$o(x)>\frac{255}{301(p_2+1)}n$, which implies that
$[G:\langle x\rangle]< \frac {301}{255}(p_2+1)<2p_2$.

First, suppose that $p_2$ does not divide $[G:\langle x\rangle]$. Then
$\langle x\rangle$
contains a cyclic Sylow $p_2$-subgroup $P$ of $G$ and as shown in the proof
of Proposition D, $P$ is a normal  cyclic Sylow subgroup of $G$.
Now Lemma 2.3(5) and our assumption imply that
$\psi(P)\psi(G/P)\geq \psi(G)=\frac {85}{301}\psi(P)\psi(|G/P|)$,
with equality if and only if $P$ is central in $G$.
Hence $\psi(G/P)\geq \frac {85}{301}\psi(|G/P|)$. Notice that $|G/P|=3m_1$,
with $m_1$ odd and $(m_1,3)=1$. 

If $G/P$ is non-cyclic, then
by Proposition D the above inequality implies that 
$\psi(G/P)= \frac {85}{301}\psi(|G/P|)$. Hence 
$\psi(P)\psi(G/P)= \psi(G)$ and by Lemma 2.3(5) $P$ is central in $G$.
Thus $G=P\times F$, where $F\simeq G/P$. Moreover, by induction 
$G/P=C_k\times(C_7\rtimes C_3)$, where $(k,42)=1$, and since 
$F\simeq G/P$, we may also write $F=C_k\times(C_7\rtimes C_3)$, 
where $(k,42)=1$. Since $(|P|,|G/P|)=1$, it follows that $(|P|,42k)=1$
and $G=C_{|P|}\times F=C_{k|P|}\times(C_7\rtimes C_3)$,
where $(k|P|,42)=1$, as required.

So assume that $G/P$ is cyclic. 
Write $G=P\rtimes F$, where $F$ is isomorphic to $G/P$.
Then $F$ is cyclic and $3$ divides $|F|$. Since $C$ is not normal in $G$, $3$
does not divide
$C_F(P)$. Therefore, by Lemma 2.3(7), $C_3$ acts fixed-point-freely on $P$,
implying that $p_2$ is congruent to $1$ (mod $3$). 
Hence $p_2\geq q_1$ and $g_3(p_2)\leq g_3(q_1)$.
Applying Lemma 2.5 we
obtain
$$\frac {\psi(G)}{\psi(|G|)}\leq \frac
{p_2^2-p_2+1+6p_2}{7(p_2^2-p_2+1)}
=g_3(p_2)\leq g_3(q_1)=\frac {\psi(G)}{\psi(|G|)}.$$
Hence $p_2=q_1=7$ and 
$\frac {\psi(G)}{\psi(|G|)}=\frac
{p_2^2-p_2+1+6p_2}{7(p_2^2-p_2+1)}$. The second equality implies,
by the proof of Lemma 2.5, 
that $|P|=p_2=7$. Since $p_2$ is the largest prime divisor of $|G|$,
it follows that $G=C_7\rtimes C_3C_{5^c}$ for some non-negative integer $c$.
Since no $5$-element can act fixed -point-freely on $C_7$, it follows by 
Lemma 2.3(7) that
$C_{5^c}$ acts trivially on $C_7$. Since $3$
does not divide $C_F(P)$, it follows that $G=C_{5^c}\times(C_7\rtimes C_3)$,
as required.

Now suppose that $p_2$ does divide  $[G:\langle x\rangle]$. Since $[G:\langle
x\rangle]<2p_2$,
it follows that $[G:\langle x\rangle]=p_2$. Then, by Theorem 3.1 of [13],
$G=P\rtimes F$, where
$P$ is a normal Sylow $p_2$-subgroup of $G$, $F$ is a cyclic subgroup of $G$,
and either
$P$ is cyclic, or $G$ is nilpotent, or
$G=(\langle a\rangle \times \langle b\rangle)\rtimes \langle y\rangle$, where
$o(a)=p_2^{\alpha -1} $ for some integer $\alpha >1$, $o(b)=p_2$,
$(o(y),p_2)=1$,
$a^y=a$ and $b^y=b^r$
for some integer $r$ not congruent to $1$ (mod $p_2$). 
Since $G/P$ is cyclic, the previous arguments imply
that if $P$ is cyclic, then $G=C_{5^c}\times(C_7\rtimes C_3)$, where $c$ is 
a non-negative integer,
as required.  Since $C$ is not normal in $G$, $G$ is not nilpotent.
Finally, since $p_2\geq 5$, in the last case we reach 
the following contradiction, using
Lemma 2.4, Theorem 9 and Lemma 2.1:
$$\psi(G)<\psi(\langle a\rangle \times \langle b\rangle)\psi(\langle
y\rangle)\leq f(p_2)\psi(|G|)\leq f(5)\psi(|G|)=\frac {121}{521}\psi(|G|)
<\frac {85}{301}\psi(|G|). $$

The proof in the case $q=3$ is now complete.
 
Now assume that $q>3$.

If $G=C_k\times C_q\times C_p\times C_p$ with $(k,p!)=1$, then
$\psi(G)=\psi(C_k)\psi(C_q)\psi(C_p\times C_p)
=\psi(C_k)\psi(C_q)(\frac {p^4+p^3-p^2+1}{p+1})$
and $\frac {\psi(G)}{\psi(|G|)}=
(\frac {p^4+p^3-p^2+1}{p+1})\frac 1{\psi(C_{p^2})}=
 \frac {p^4+p^3-p^2+1}{p^5+1}$, as claimed.

Now we turn to the "only if" part.
So suppose that $q>3$, $|G|=n=qm$, $m$ is odd, $(m,q)=1$ and     
$\frac {\psi(G)}{\psi(|G|)}=\frac {p^4+p^3-p^2+1}{p^5+1}=f(p)$.
Then $G$ is non-cyclic and $G=M\rtimes C$, where $C\simeq C_q$ 
and $q$ does not divide $|M|$.  Denote by $r$ the least prime divisor 
of $|M|$. Clearly $r\geq p$.
Our aim is to show that $G=C_k\times C_q\times C_p\times C_p$ with $(k,p!)=1$. 
We argue by induction on $|G|$.

If $G=M\times C$, then $M$ is non-cyclic and $\psi(G)=\psi(M)\psi(C)$. 
Moreover, since $r\geq p$, Theorem 9 and Lemma 2.1 imply that 
$\psi(M)\leq f(r)\psi(|M|)\leq f(p)\psi(|M|)$.
Hence 
$\frac {\psi(G)}{\psi(|G|)}=\frac {\psi(M)\psi(C)}{\psi(|M|)\psi(C)}
\leq f(p)$. But by our
assumptions 
$\frac {\psi(G)}{\psi(|G|)}=f(p)$, so $r=p$ and
$\frac {\psi(M)}{\psi(|M|)}=f(p)=f(r)$. Therefore, by Theorem 9,
$M=C_p\times C_p\times C_{k}$, with $(k,p!)=1$ and hence 
$G=C_q\times C_p\times C_p\times C_k$ with $(k,p!)=1$, as required.

Suppose, now, that $C$ is not normal in $G$.

Let $p_2$ be the largest prime divisor of $|G|$ and let $P$ be a Sylow 
$p_2$-subgroup of $G$. Then, by Lemma 2.3(4), we have
$\psi(|G|)\geq \frac q{p_2+1}n^2$ and our assumptions imply that
$\psi(G)=f(p)\psi(|G|)\geq \frac{f(p)q}{p_2+1}n^2$.
Hence there exists $x\in G$ such that
$o(x) >\frac{f(p)q}{p_2+1}n$, which implies that
$[G:\langle x\rangle]< \frac {p_2+1}{f(p)q}$.
Now by the Bertrand's conjecture we have $p\leq 2q-3$, so
$q\geq \frac {p+3}2$. Since 
$f(p)=\frac {p^4+p^3-p^2+1}{p^5+1}>\frac 1p$, it follows that
$f(p)q>\frac {p+3}{2p}$ and as $p_2\geq p$, we get
$$[G:\langle x\rangle]< \frac {2p}{p+3}(p_2+1)<2p_2.$$ 
 
First, suppose that $p_2$ does not divide $[G:\langle x\rangle]$. Then
$\langle x\rangle$
contains a cyclic Sylow $p_2$-subgroup $P$ of $G$. As shown in the proof
of Proposition D, $P$ is a normal  cyclic Sylow subgroup of $G$.
Now Lemma 2.3(5) and our assumptions imply that 
$\psi(P)\psi(G/P)\geq \psi(G)=
f(p)\psi(P)\psi(|G/P|)$.
Hence $\psi(G/P)\geq f(p)\psi(|G/P|)$. Notice that $|G/P|=qm_1$,
with $m_1$ odd and $(m_1,q)=1$.

Suppose, first, that $G/P$ is non-cyclic.
Then, by Proposition D, the above inequality implies that
$\psi(G/P)= f(p)\psi(|G/P|)$. Hence
$\psi(P)\psi(G/P)=\psi(G)$ and by Lemma 2.3(5) $P$ is central in $G$.
Thus $G=P\times F$, where $F\simeq G/P$. Moreover, by induction,
$G/P=C_q\times C_p\times C_p\times C_{k_1}$ with $(k_1,p!)=1$ and since
$F\simeq G/P$, we may also write $F=C_q\times C_p\times C_p\times C_{k_1}$
with $(k_1,p!)=1$. As $(|P|,|G/P|)=1$, we have $(|P|,k_1pq)=1$.
Hence $(|P|,k_1)=1$ and as $p_2>p$, also $(|P|,p!)=1$.  
Thus $G=P\times F=C_{k_1|P|}\times C_q\times C_p\times C_p$
with $(k_1|P|,p!)=1$, as required.

Suppose, next, that $G/P$ is cyclic.
Write $G=P\rtimes F$, where $F$ is isomorphic to $G/P$.
Then $F$ is cyclic and $q$ divides $|F|$.
Since $C$ is not normal in $G$, $q$
does not divide
$C_F(P)$. Therefore, by Lemma 2.3(7), $C_q$ acts fixed-point-freely on $P$,
implying that $p_2$ is congruent to $1$ (mod $q$).
Hence  $p_2\geq p_1$ and $g_q(p_2)\leq g_q(p_1)$.
Applying Lemma 2.5 we
obtain
$$\frac {\psi(G)}{\psi(|G|)}\leq \frac
{p_2^2-p_2+1+p_2(q^2-q)}{(q^2-q+1)(p_2^2-p_2+1)}
=g_q(p_2)\leq g_q(p_1).$$
By Proposition 2.2(1) $g_q(p_1)<f(p)$, so
$\frac {\psi(G)}{\psi(|G|)}<f(p)$, in contradiction to our assumptions.

Now suppose that $p_2$ does divide  $[G:\langle x\rangle]$. Since $[G:\langle
x\rangle]<2p_2$,
it follows that $[G:\langle x\rangle]=p_2$. Then, by Theorem 3.1 of [13],
$G=P\rtimes F$, where
$P$ is a normal Sylow $p_2$-subgroup of $G$, $F$
is a cyclic subgroup of $G$,
and either
$P$ is cyclic, or $G$ is nilpotent, or
$G=(\langle a\rangle \times \langle b\rangle)\rtimes \langle y\rangle$, where
$o(a)=p_2^{\alpha -1} $ for some integer $\alpha >1$, $o(b)=p_2$,
$(o(y),p_2)=1$,
$a^y=a$ and $b^y=b^r$
for some integer $r$ not congruent to $1$ (mod $p_2$). 

If $P$ is cyclic, 
then $G=P\rtimes F$, where both $P$ and $F$ are cyclic, and we reach a 
contradiction as above.
Since $C$ is not normal in $G$, $G$ is not nilpotent.
Finally, consider the last case. Since $p_2\geq p$, we have 
$f(p_2)\leq f(p)$. Now, applying 
Lemma 2.4 and Theorem 9, we reach the following
contradiction:
$$\psi(G)<\psi(\langle a\rangle \times \langle b\rangle)\psi(\langle
y\rangle)\leq f(p_2)\psi(|G|)\leq f(p)\psi(|G|)=\psi(G).$$

The proof of Proposition E is now complete.
\qed
\enddemo

\heading 5. Proof of Corollary B \\
\endheading
\demo{Proof} Recall that $f(x)=\frac {x^4+x^3-x^2+1}{x^5+1}$.
Let $G$ be a non-cyclic group of odd order $n=qm$, 
where $q$ is the minimum
prime divisor of $n$ and $(m,q)=1$.

Suppose, first, that $q>3$. Then, by Theorem A,
$$\frac {\psi(G)}{\psi(|G|)}\leq f(p),$$
where $p$ is the smallest prime bigger than $q$. 
Since $q>3$, it follows that $q\geq 5$, $p\geq 7$ and by Lemma 2.1
$f(p)\leq f(7)=\frac {337}{2101}$.
Hence if $q>3$, then
$$\frac {\psi(G)}{\psi(|G|)}\leq f(7)=\frac {337}{2101},$$
and by Proposition E equality holds if and only if $n=5\cdot 7^2\cdot m_1$ 
with
$(m_1,7!)=1$ and $G=C_5\times C_7\times C_7\times C_{m_1}$,
as required.

Suppose, now, that $q=3$. Then by Theorem A
$$\frac {\psi(G)}{\psi(|G|)}\leq \frac {85}{301},$$
with equality if and only if $n=3\cdot 7\cdot m_1$ with $(m_1,42)=1$
and $G=(C_7\rtimes C_3)\times C_{m_1}$. 

Since $\frac{337}{2101} < \frac{85}{301}$, the result concerning groups with $q \geq 3$ follows immediately.

\qed
\enddemo

\enddocument